\documentclass[11pt]{amsart}
\usepackage{geometry}                
\usepackage{graphicx}
\usepackage{amssymb}
\usepackage{epstopdf}
\usepackage{color}
\usepackage{tikz}
\usetikzlibrary{arrows, decorations.markings, decorations.pathmorphing, backgrounds, positioning, fit, petri}

\newtheorem{thm}{Theorem}[section]
\newtheorem{cor}[thm]{Corollary} 
\newtheorem{lemma}[thm]{Lemma}
\newtheorem{example}[thm]{Example}
\newtheorem{defn}[thm]{Definition}  
\newtheorem{rmk}[thm]{Remark}

\title{Directed diagrammatic reducibility}
\author{Jens Harlander and Stephan Rosebrock}
\date{\today}                                           

\begin{document}
\begin{abstract}
We introduce the notion of directed diagrammatic reducibility which is a relative version of diagrammatic reducibility. Directed diagrammatic reducibility has strong group theoretic and topological consequences. A multi-relator version of the Freiheitssatz in the presence of directed diagrammatic reducibility is given. Results concerning asphericity and $\pi_1$-injectivity of subcomplexes are shown. We generalize the Corson-Trace characterization of diagrammatic reducibility to directed diagrammatic reducibility. We compare diagrammatic reducibility of relative presentations to directed diagrammatic reducibility. Classical tools for showing diagrammatic reducibility, such as the weight test, the max/min test, and small cancellation techniques are adapted to directed diagrammatic reducibility. The paper ends with some applications to labeled oriented trees. \end{abstract}
\maketitle


\section{Introduction} A 2-complex $K$ is called {\it aspherical} if $\pi_2(K)=0$.
There are combinatorial versions of asphericity, the strongest being diagrammatic reducibility, or DR for short. Recall that a map between 2-complexes is called {\em combinatorial} if it maps open cells homeomorphically to open cells. For a 2-complex $K$, a combinatorial map $f\colon C \to K$ is called a {\it spherical diagram}, if $C$ is a 2-sphere with cell structure. A 2-complex $K$ is {\it diagrammatically reducible}, DR for short, if every spherical diagram $f\colon C \to K$ contains an edge $e$ so that the 2-cells in $C$ that contain $e$ in their boundary map to the same 2-cell in $K$ with opposite orientations by folding over $e$. The edge $e$ is called a {\it folding edge} in $C$. Diagrammatic reducibility is a 2-dimensional version of free reductions of cycles in a graph. It was first considered by Sieradski \cite{Sie83}. See Gersten \cite{Ger87} for a good overview of diagrammatic techniques with applications to combinatorial topology and group theory. \\

A presentation $P$ defines a 2-complex $K(P)$ and a group $G(P)$ which is the fundamental group of $K(P)$ by Van Kampen's theorem.
A presentation $P$ is called DR if its associated 2-complex $K(P)$  is DR.
For a set $X$ call a subset $S$ {\it proper} if $S\ne X$ ($S$ may be empty). If $S$ is a proper subset of the set of generators of $P$ let $P_S$ be the sub-presentation of $P$ carried by $S$. We say {\em $P$ is DR directed away from $S$} if every spherical diagram $f\colon C\to K(P)$ that is not already a diagram over $K(P_S)$ contains a folding edge that is mapped to an edge in $K(P)$ that is not an edge in $K(P_S)$. This is a relative version of diagrammatic reducibility closely related to (but different from) relative vertex asphericity which was recently developed by the authors in \cite{HaRo18}.\\

Here is an outline of the paper. In Section \ref{sec:basics} we show that if $P$ is DR directed away from $S$ then $\pi_2(K(P))$ is generated (as a $G(P)$-module) by the image of the inclusion induced map $\pi_2(K(P_S))\to \pi_2(K(P))$ and the inclusion $K(P_S)\to K(P)$ is $\pi_1$-injective. We provide basic examples and give, as our first application, a multi-relator Freiheitssatz for presentations $P$ that are DR in all directions. In Section \ref{sec:corsontrace} we characterize directed diagrammatic reducibility in terms of finite subcomplexes of $\tilde K(P)$, the universal covering of $K(P)$. In the standard DR setting this was done by Corson and Trace \cite{CT00}. In Section \ref{sec:relpres} diagrammatic reducibility of relative presentations introduced by Bogley and Pride in \cite{BP92} is compared with directed diagrammatic reducibility. We show that the two concepts are equivalent. However, the viewpoints are very different. In the Bogley/Pride setting, the starting point is a relative presentation where the generators are a free product $H*F({\bf x})$, where $H$ is fixed. We start with an ordinary presentation on generators ${\bf x}$ and consider DR directed away from a proper subset $S\subseteq {\bf x}$, where $S$ is not necessarily fixed but can vary over ${\bf x}$. This enables us to consider the idea of diagrammatic reducibility in all directions. We can formulate a natural multi-relator Freiheitssatz in Corollary \ref{cor:freiheit}. Retaining combinatorial information about the subpresentation $P_S$ carried by the subset $S$ of generators of $P$ rather than only keeping its fundamental group and passing to a relative presentation, might give an advantage in designing tools for detecting DR in the two settings. Such tools, a max/min and a weight test for showing directed diagrammatic reducibility are given in Section \ref{sec:tests}.\ Theorem \ref{s44} and the examples we exhibit at the end of Section \ref{sec:appl} might have been overlooked from a purely relative presentation point of view.\\

We conclude this introduction with some remarks about diagrams. In this paper we consider surface diagrams where the surface can be different from a sphere. A {\em surface diagram over a 2-complex $K$} is a combinatorial map $f\colon F\to K$, where $F$ is a compact surface, possibly with boundary, with a cell-structure.  If $c$ is a cell in $F$, we call $f(c)$ the {\em label} of $c$. The labeled cell-complex $F$ contains all the information of the combinatorial map $f$, and we sometimes refer to the labeled $F$ as a surface diagram over $K$. A {\em Van Kampen diagram over K} is a combinatorial map $M\to K$ where $M$ is a simply-connected planar region with a cell structure. Topologically, $M$ is a ``tree of discs". It is well known that an edge loop $\gamma$ in $K$ represents the trivial element in the fundamental group if and only there exists a reduced Van Kampen diagram $M\to K$ whose boundary loop is mapped to $\gamma$. 

We are primarily concerned with 2-complexes and groups defined by a finite presentation 
\[ P=\langle x_1,\ldots ,x_n\mid r_1,\ldots ,r_m\rangle.\]
Recall that $K(P)$ is the 2-complex with one vertex, edges and 2-cells in one-to-one correspondence  with the generators and relators respectively. The 2-cell $\Delta(r)$ for the relator $r$ is attached to the one skeleton according to the word $r$. The universal cover $\tilde K(P)$ has as vertex set the elements of $G(P)$, oriented edges $(g,x_i)$, $g\in G(P)$, $1\le i \le n$, and 2-cells $(g,r_j)$, $g\in G(P)$, $1\le j \le m$. The 1-skeleton of $\tilde K(P)$ is also referred to as the Cayley graph of $G(P)$ associated with the set of generators. The boundary of $(g,r_j)$ is the lift of the path $r_j$ to the Cayley graph, starting (and ending) at $g$.

\vspace{5ex}

\section{Directed diagrammatic reducibility}\label{sec:basics} 

As was mentioned in the introduction, a 2-complex $K$ is called diagrammatically reducible (or DR for short) if each spherical diagram $f\colon C\to K$ contains a folding edge.

\begin{defn}\label{defDrdira} Let $P=\langle x_1,\ldots ,x_n\mid r_1,\ldots ,r_m\rangle$ be a presentation and $S$  a  proper subset (possibly empty) of the set of generators. We say that $P$ is
\begin{itemize}
\item {\em DR directed away from $S$} if every spherical diagram $f\colon C\to K(P)$ that contains an edge with label not from $S$  also contains a folding edge with label not from $S$;
\item {\em DR in all directions} if every spherical diagram $f\colon C\to K(P)$ that contains an edge labeled $x_i$ also contains a folding edge with label $x_i$, $i=1,\ldots ,n$. Note that this implies that $P$ is DR directed away from all proper subsets $S$.
\end{itemize}
\end{defn}

\noindent If $S=\emptyset$ then DR directed away from $S$ simply means DR.\\[2ex]

\noindent{\bf Examples and Comments:}\\

\noindent 1) Just like DR, directed DR is a hereditary property: Let $T$ be a sub-presentation of a presentation $P$ and let $S$ be a subset of the set of generators of $T$. Then if $P$ is DR directed away from $S$ so is $T$. If $P$ is DR in all directions, then so is $T$.

However, a presentation $P$ can be DR directed away from $S$ but not DR directed away from some $S_0\subset S$. Here is an example. Suppose that 
$$P=\langle a, b, c\mid aba^{-1}b^{-2}, bab^{-1}a^{-2}, cab\rangle.$$ 
Note that $G(P)$ is the trivial group. $P$ is DR away from $S=\{a,b\}$. This is because if a spherical diagram $f\colon C\to K(P)$ contains an edge $e$ with label not in $S$, then the label on $e$ is $c$. Since $c$ is a free edge in $K(P)$, this implies that $e$ is a folding edge in $C$. 
Note that $P$ is not DR directed away from $S_0=\{ a \}$ by Theorem \ref{RelPiS} (2) below. The presentation is also not DR because the  subpresentation $P_0=\langle a, b \mid  aba^{-1}b^{-2}, bab^{-1}a^{-2}\rangle$ is not DR. The 2-complex $K(P_0)$ is simply connected but does not have a free edge (see Corson-Trace \cite{CT00}, and also Section \ref{sec:corsontrace} in this article).\\ 

\noindent 2) Subdivisions preserve DR but not directed DR. Consider $P=\langle a, b \mid aba^{-1}b^{-1}\rangle$ and $P'=\langle a, b, c \mid abc, c^{-1}a^{-1}b^{-1}\rangle$. Then $P$ is DR in all directions (see the next example) but $P'$ is not DR directed away from $S=\{ a, b\}$: the two triangles coming from the relators can be glued together to form a disc-diagram that contains an edge labeled $c$ but not a folding edge labeled $c$. Thus $P'$ is not DR directed away from $S$ by Theorem \ref{RelPiS} (2). \\

\noindent 3) Direct products of graphs are DR in all directions. Consider
$$P=\langle x_1,\ldots ,x_m, y_1,\ldots ,y_n \mid [x_i,y_j], 1\le i\le m, 1\le j\le n \rangle.$$ 
The presentation $P$ is DR in all directions. Let $f\colon C\to K(P)$ be a spherical diagram that contains an edge labeled by $y_k$. We will show that there is a folding edge labeled by $y_k$. If $d$ is a 2-cell in $C$ with two of its edge labels $y_k$, then draw a red line across $d$ connecting the midpoints of these two edges. The red graph in $C$ we draw in this way consists of red circles. Take one of these red circles and consider the gallery of 2-cells that comes with it. This gallery forms an annulus $A$, the two boundary circles are labeled by a word $w$ in the $x$'s. Edges connecting the two boundary circles are all labeled by $y_k$. Since the subgroup of $G(P)$ generated by $x_1,\ldots ,x_m$ is free, the word $w$ must contain a cancelling pair $x_i^{\epsilon}x_i^{-\epsilon}$. Thus the annulus $A$ contains a cancelling pair of 2-cells labeled by $[x_i,y_k]^{\epsilon} [x_i,y_k]^{-\epsilon}$ with a folding edge labeled by $y_k$. The argument for an $x_k$ is analogues.\\

\noindent 4) Orientable closed surfaces are DR in all directions. Let 
$$P=\langle x_1, x_2,\ldots ,x_{2g-1}, x_{2g}\mid [x_1,x_2]\ldots [x_{2g-1}, x_{2g}]\rangle$$ 
be the standard presentation of the fundamental group of an orientable surface of genus $g$. Then $P$ is DR in all directions. Let $f\colon C\to K(P)$ be a spherical diagram that contains an edge labeled by $x_1$, say. We will show that there is a folding edge labeled by $x_1$. Just as in the previous example, $C$ contains an annulus $A$ with the two boundary components labeled by words $u$ and $v$ in $\{ x_2,\ldots ,x_{2g}\}$. Edges connecting the boundary components are all labeled with $x_1$. Since the subgroup of $G(P)$ generated by $x_2,\ldots ,x_{2g}$ is free, the words $u$ and $w$ must contain a cancelling pair of letters. Hence the annulus $A$ contains a cancelling pair of 2-cells with folding edge labeled by $x_1$. \\

Let $P=\langle x_1,\ldots ,x_n\mid r_1,\ldots ,r_m\rangle$ be a presentation and $S$ a subset of the generators.

\begin{thm}\label{RelPiS} Suppose that $P$ is DR directed away from the proper subset $S$ of the generators. Then 
\begin{enumerate}
\item $\pi_2(K(P))$ is generated (as a $G(P)$-module) by the image of the inclusion induced map $\pi_2(K(P_S))\to \pi_2(K(P))$; furthermore
\item every disc diagram $g\colon D\to K(P)$ with boundary labeled by a word in $S$, that contains a label not from $S$, has a folding edge with label not from $S$. Consequently, the inclusion induced map $\pi_1(K(P_S))\to \pi_1(K(P))$ is injective.
\end{enumerate}
\end{thm}

\noindent Proof. Suppose $f\colon C\to K(P)$ is a reduced spherical diagram. If $f(C)$ is not contained in $K(P_S)$ then $C$ contains an edge $e$ so that $f(e)\not\in S$. Since we assumed that $P$ is DR directed away from $S$ it follows that $C$ contains a folding edge $e'$ so that $f(e')\not\in S$, contradicting the assumption that $f\colon C\to K(P)$ is reduced. Since $\pi_2(K(P))$ is generated (as a $G(P)$-module) by reduced spherical diagrams, the first statement follows.

Suppose $g\colon D\to K(P)$ is a disc diagram as in statement (2). We double $D$ and construct a spherical diagram $g'\colon C=D_1\cup D_2 \to K(P)$, where $D_1$ is mapped by $g$ and $D_2$ is mapped by $-g$ (an orientation reversion followed by $g$).  Note that $C$ contains an edge with label not in $S$. Since $P$ is DR away from $S$ this spherical diagram contains a folding edge with label not from $S$. This folding edge can not occur on $\partial D_1=\partial D_2$. Thus $D_1$ or $D_2$ contain an interior folding edge with label not in $S$ and hence so does $D$.

We next show $\pi_1$-injectivity. Suppose $w$ is a word in $S^{\pm 1}$ that represents a non-trivial element of $\pi_1(K(P_S))$ that maps to a trivial element in $\pi_1(K(P))$. Then there exists a reduced Van Kampen diagram $f\colon M\to K(P)$ where the boundary of $M$ is labeled by $w$ and is mapped to $K(P_S)$. Here $M$ is a planar simply connected region with a cell structure. Note that $M$ is a tree with discs attached at some vertices. One of these discs, say $\bar D$, contains an edge $e$ such that $f(e)\not\in K(P_S)$, otherwise $f(M)\subseteq K(P_S)$ which would imply that $w=1$ in $\pi_1(K(P_S))$. Thus ${\bar g}=f|_{\bar D}\colon {\bar D}\to K(P)$ is a disc diagram as in statement (2) and hence contains a folding edge, contradicting the fact that we assumed $f\colon M\to K(P)$ is reduced. \qed

\bigskip\noindent We recall the Freiheitssatz for 1-relator groups: Suppose $P=\langle x_1,\ldots ,x_n\mid r \rangle$ is a 1-relator presentation, where $r$ is a cyclically reduced word that contains all the generators. Then any proper subset $S$ of $\{x_1,\ldots ,x_n\}$ generates a free subgroup of $G(P)$ with basis $S$. The following three results should be viewed as multi-relator versions of this celebrated result.

\begin{thm} Let $P=\langle x_1,\ldots ,x_n\mid r_1,\ldots ,r_m\rangle$ be a presentation that is DR in all directions. Then the inclusion induced homomorphism $G(P_S)\to G(P)$ is injective for every subset $S$ of the generators. \end{thm}

\noindent Proof. If $S$ is the set of generators of $P$ then $G(P_S)=G(P)$ and the statement is true. If $S$ is a proper subset of the set of generators then the result follows from Theorem \ref{RelPiS} using the fact that $P$ is DR directed away from $S$.\qed

\begin{thm} Let $P=\langle x_1,\ldots ,x_n\mid r_1,\ldots ,r_m\rangle$ and $S$ a proper subset of the generators. Assume that each $r_i$ contains a generator not from $S$. If $P$ is DR
directed away from $S$, then $S$ generates a free subgroup of $G(P)$ with basis $S$.
\end{thm}

\noindent Proof. Since each $r_i$ contains a generator not from $S$ we have that $P_S=\langle S\mid\ \rangle $ and $G(P_S)$ is free. Now Theorem \ref{RelPiS} (2) gives the desired result.\qed\\

\begin{cor}\label{cor:freiheit} Let $P=\langle x_1,\ldots ,x_n\mid r_1,\ldots ,r_m\rangle$ be a presentation where each $r_i$ contains all the generators. If $P$ is DR in all directions, then any proper subset $S$ of $\{x_1,\ldots ,x_n\}$ generates a free subgroup of $G(P)$ with basis $S$. 
\end{cor}

We end this section with a strengthening of the classical Freiheitssatz for torsion-free 1-relator groups. Recall that a group $G$ is {\it left orderable} if there exists an order relation $<$ so that for all $x\ne y$ in $G$ either $x<y$ or $y<x$. Furthermore, if $x<y$ and $z\in G$, then $zx<zy$. One relator groups are locally indicable  by a result of Brodskii \cite{Brodskii} (see also Howie \cite{Howie}). A theorem of Burns-Hale \cite{BurnsHale} implies that locally indicable groups are left orderable. Our use of orderability in the proof of the next result is related to key ideas in Barmak-Minian \cite{BM18}.

\begin{thm}\label{thm:1rel}Let $P=\langle x_1,\ldots ,x_n \mid r \rangle$ be a one-relator presentation of a group $G$ where $r$ is a cyclically reduced relator that is not a proper power. Then $P$ is DR in all directions.\end{thm}

\noindent Proof: Choose a left order on $G(P)$. $K=K(P)$ is a 2-complex with a single 2-cell $d$. Let $\tilde d$ be a fixed lift of $d$ in the universal cover $\tilde K$. Note that if $r$ has length $k$, then $\tilde d$ is attached by a simple edge loop $\tilde r$ of length $k$. Otherwise $r$ would contain a non-trivial proper subword that represents the trivial element in $G(P)$. This is not possible by a result of Weinbaum (see Lyndon-Schupp \cite{LS77}, Proposition 5.2.9, page 110). Thus all $k$ edges that occur in the attaching path $\tilde r$ have distinct initial vertices. Define
$$m(x_i)=\max\{ g\in G\mid (g,x_i)\ \mbox{is an edge in the attaching loop $\tilde r$ of $\tilde d$}\}.$$
It follows from the remarks above that for every $x_i$ that occurs in $r$ we have a unique edge $(m(x_i), x_i)$ in the boundary $\tilde r$ of $\tilde d$.\\

Suppose $f\colon C\to K$ is a spherical diagram that contains an edge labeled $x_j$.  Consider a lift $\tilde f\colon C\to \tilde K$. Let $(h,x_j)$ be an edge in $\tilde f(C)$ so that $h$ is maximal among all the $(g,x_j)$ edges in $\tilde f(C)$. Let $g_1\tilde d$ and $g_2\tilde d$ be the 2-cells in $\tilde f(C)$ that contain $(h,x_j)$. Then $(g_1^{-1}h,x_j)=(m(x_j),x_j)$ and $(g_2^{-1}h,x_j)=(m(x_j),x_j)$. Therefore $g_1=g_2$. It follows that the edge in $C$ that maps to $(h,x_j)$ is a folding edge. \qed \\

Thus the Freiheitssatz for torsion-free 1-relator groups ``follows" from Corollary \ref{cor:freiheit}. However, the proof of Theorem \ref{thm:1rel} uses a deep result of Brodskii (local indicability of torsion-free 1-relator groups), and Brodskii's proof relies on the Freiheitssatz.

\vspace{5ex}

\section{Corson-Trace for directed diagrammatic reducibility}\label{sec:corsontrace}

If $K(P)$ is the 2-complex build from the presentation $P$ let $\tilde K(P)$ be its universal cover and $\tilde K(P)^{(1)}$ the 1-skeleton of the universal cover. In \cite{CT00} Corson and Trace have shown the following result:

\begin{thm}\label{thm:classicCT} A presentation $P$ is DR if and only if every finite subcomplex of $\tilde K(P)$ collapses into the 1-skeleton $\tilde K(P)^{(1)}$.
\end{thm}

\noindent In this section we retell the Corson-Trace story in the context of directed DR. 

\begin{lemma}\label{RedSurDiaS} Let $P$ be a finite presentation and let $S_0$ be the subset of the generators which correspond to the free edges in $K(P)$. Let $x$ be a generator contained in a relator but not in $S_0$. Then there exists a closed oriented surface diagram $f\colon F\to K(P)$ that contains an edge labeled by $x$, which is reducible only at edges labeled by elements from $S_0$.
\end{lemma}

\noindent Proof. To avoid an excessive amount of subscripts we use letters $a, b, c,\ldots $ for the set of generators in $P$. For every relator $r$ draw two $|r|$-gons $\Delta(r)$ and $\Delta(r^{-1})$ in the plane. Orient the edges of these $|r|$-gons and label them according to the relator $r$, using subscripts. For example if $r=abab^{-1}ac$ then label the edges of $\Delta(r)$ clockwise by 
$$a_1(r),b_1(r),a_2(r),b_2^{-1}(r),a_3(r),c_1(r)$$ and $\Delta(r^{-1})$ counter-clockwise by 
$$a_1(r^{-1}),b_1(r^{-1}),a_2(r^{-1}),b_2^{-1}(r^{-1}),a_3(r^{-1}),c_1(r^{-1}).$$
Note that the number of edges on the $2m$ discs $\Delta(r_i)$, $\Delta(r_i^{-1})$, $i=1,...,m$,  with labels of the form $a_{*}(*)$
is even, equal to 2 if $a$ is a free edge of $K(P)$, otherwise greater or equal to $4$. The same holds true for every other generator letter. We now match edges of the $2m$ discs in pairs. If $a\in S_0$ is a free edge that occurs only once among the set of relators, say in $r$, then match $a_1(r)$ and $a_1(r^{-1})$. If $a$ is not a free edge then match $a_i(r_k)$ with any $a_j(r_l^{\pm 1})$, but avoid the match with $a_i(r_k^{-1})$. One concrete way to do this is the following: Define
\begin{eqnarray*}
L_+(a) & : & \ a_1(r_1), a_2(r_1),\ldots ,a_1(r_2), a_2(r_2),\ldots ,a_1(r_3),a_2(r_3),\ldots\\
L_{-}(a) & : & \ a_1(r_1^{-1}), a_2(r_1^{-1}),\ldots ,a_1(r_2^{-1}), a_2(r_2^{-1}),\ldots ,a_1(r_3^{-1}),a_2(r_3^{-1}),\ldots
\end{eqnarray*}
Cyclically match item $N$ on the $L_+(a)$ list with item $N+1$ on the $L_-(a)$ list. Proceed in this fashion with all generators $a, b, c,\ldots$.
This way we build a surface, possibly consisting of more than one component. Let $F_0$ be a component that contains an edge $x_i(r_k)$. If we erase from the edge labels  all information except the letter (for example $a_i(r_j^{\pm 1})$ becomes $a$), we obtain a closed surface diagram $f_0\colon F_0\to K(P)$ that contains the label $x$. Passing to the 2-fold orientable covering $F$ of $F_0$ we obtain a closed orientable surface diagram $f\colon F\to K(P)$ that contains the label $x$. By construction $F$ is reducible only at edges with labels from $S_0$. \qed \\

Here is a topological version of this lemma.

\begin{lemma}\label{TopRedSurDiaS}
Let $X$ be a 2-complex and let $E_0$ be the set of boundary edges of $X$. Let $e$ be a fixed edge that is part of a 2-cell of $X$ and is not a boundary edge. Then there exists a closed oriented surface diagram $F\to X$ so that $F$ contains an edge labeled by $e$ and $F$ is reducible only at edges which carry labels from $E_0$.
\end{lemma}

We recall the well known fact that if $P$ is a presentation and $w$ a word in the generators of $P$ that represents the trivial element in $G(P)$, then there exists a reduced Van Kampen diagram $M_{w}\to K(P)$ ($M_w$ is a simply-connected planar region) with boundary word $w$. Here is the Corson-Trace result for directed DR. The original Theorem \ref{thm:classicCT} is covered by the case $S=\emptyset$.

\begin{thm}\label{thm:GenCT}The presentation $P$ is DR directed away from $S$ if and only if every finite subcomplex $X \subseteq \tilde K(P)$
can be collapsed into $p^{-1}(K(P_S))\cup \tilde K(P)^{(1)}$, where 2-cells are collapsed across edges not in $p^{-1}(K(P_S))$ (these are edges of the form $(g,x)$, $x\notin S$). \end{thm}

\noindent Proof. Suppose first that $P$ is DR directed away from $S$. Let $X$ be a finite subcomplex of $\tilde K(P)$ that is not already contained in $p^{-1}(K(P_S))\cup \tilde K(P)^{(1)}$. We may assume every edge of $X$ is part of a 2-cell. It suffices to show that $X$ has a boundary edge that is not contained in $p^{-1}(K(P_S))$. We will show this by contradiction. Let $\tilde E_0$ be the set of boundary edges of $X$ and assume that $\tilde E_0\subseteq p^{-1}(K(P_S))$. Then $p(\tilde E_0)=S_0\subseteq S$.
By Lemma \ref{TopRedSurDiaS} there exists a closed oriented surface diagram $\tilde f\colon F\to X\subseteq \tilde K(P)$ so that $\tilde f(F)$ is not contained in $p^{-1}(K(P_S))$ and folding edges in $F$ (should there be any) are mapped to $E_0$.
We consider the composition $$f\colon F\to X\subseteq \tilde K(P) \stackrel{p}\to K(P).$$ Note that $f\colon F\to K(P)$ is a closed oriented surface diagram so that $f(F)$ is not contained in $K(P_S)$ and folding edges in $F$ (should there be any) are labeled by elements from $S_0\subseteq S$. Since $\tilde f$ is a lift of $f$, closed edge paths in $F$ are labeled by words $w$ in the generators of $P$ that represent the trivial element in $G(P)$. 
We now proceed as in the proof of Lemma 2.1 of Corson-Trace \cite{CT00}. Choose a complete set of closed cutting curves $\gamma_1,\ldots,\gamma_k$ in $F$ (cutting along all $\gamma_i$ would transform $F$ into a simply-connected region of the plane). Since the $f(\gamma_i)$ are trivial in $G(P)$ we have a Van Kampen diagram $M_{\gamma_i}$ for each $\gamma_i$. Form a 2-complex $L$ by attaching $M_{\gamma_i}$ to $F$ for every $\gamma_i$. We have
\begin{enumerate}
\item $L\subseteq S^3$ is a simply connected 2-skeleton of a cell decomposition of the 3-sphere $S^3$;
\item we have combinatorial maps $F\stackrel{\alpha}{\rightarrow}L\stackrel{\beta}{\rightarrow}K(P)$ so that $f=\beta\circ \alpha$.
\end{enumerate}
Let $L'$ be a 2-complex that has the minimal number of 2-cells among all 2-complexes that satisfy the two conditions just stated and let $\beta'\colon L'\to K(P)$ be the corresponding combinatorial map. Note that $L'$ can not contain a 2-cell $d'$ with a free edge $e'$ with label $x\notin S$. Suppose it does. Then there is no edge $e$ in $F$ so that $\alpha(e)=e'$, because such an edge $e$ would be a folding edge in $F$ with label $x\notin S$, which we know does not exist. Thus we can collapse $d'$ in $L'$ and contain a 2-complex that satisfies our two conditions and contains fewer 2-cells than $L'$. A contradiction to minimality. Note that we do not rule out free edges in $L'$ with label in $S$.

Consider a spherical diagram $h'\colon C\to L'$ that arises as an attaching map of a 3-cell of $S^3$ so that the composition $h=\beta'\circ h'\colon C\to K(P)$ does not entirely map into $K(P_S)$. Such an attaching map exists because otherwise $\beta'(L')$ and hence $f(F)$ would be entirely contained in $K(P_S)$, which is not the case. We will show that the spherical diagram $h\colon C\to K(P)$ can not contain a folding edge with label $x\notin S$. Suppose it does. Let $e'$ be an edge in $C$ with label $x\notin S$, and let $d'_1$ and $d'_2$ be the folding pair of 2-cells for $e'$. Both $d'_1$ and $d'_2$ map to the same 2-cell $d$ in $K(P)$ under $h$. We can identify $d'_1$ and $d'_2$ to a single 2-cell $d'$ (fold $d'_1$ onto $d'_2$) and obtain a new $L'$ that satisfies our conditions and has one less 2-cell. A contradiction to minimality. The existence of the diagram $h\colon C\to K(P)$ contradicts the hypothesis that $P$ is DR directed away from $S$. So the assumption $p(E_0)=S_0\subseteq S$ is false.

Now we prove the other direction. Let $f\colon C\to K(P)$ be a spherical diagram so that $f(C)$ is not contained in $K(P_S)$. Let $\tilde f\colon C\to \tilde K(P)$ be a lift and let $X=\tilde f(C)$. Then $X$ is not contained in $p^{-1}(K(P_S))\cup \tilde K(P)^{(1)}$ and hence contains a boundary edge $\tilde e$ not in $p^{-1}(K(P_S))$. The edge $e$ in $C$ such that $\tilde f(e)=\tilde e$ is a folding edge in $f\colon C\to K(P)$ and $f(e)=x\notin S$. \qed

\begin{cor} Suppose $P$ presents a finite group. Then $P$ is $DR$ directed away from $S$ if and only if $K(P)$ collapses into $K(P_S)$, where 2-cells are collapsed across edges corresponding to generators $x\notin S$.
\end{cor}

\noindent Proof: Since we assume that $G(P)$ is finite $\tilde K(P)$ is a finite complex. By Theorem \ref{thm:GenCT} $P$ is DR directed away from $S$ if and only if $\tilde K(P)$ collapses into $p^{-1}(K(P_S))\cup \tilde K(P)^{(1)}$, where 2-cells are collapsed across edges not in $p^{-1}(K(P_S))$, i.e. edges of the form $(g,x)$, $x\notin S$. A collapsing sequence in the universal cover gives a collapsing sequence downstairs that collapses $K(P)$ into $K(P_S)\cup K(P)^{(1)}$, where 2-cells are collapsed across edges corresponding to generators $x\notin S$. Since $K(P)^{(1)}$ is a wedge of circles and $G(P)$ is assumed to be finite, we have collapsed $K(P)$ into $K(P_S)$. \qed

\vspace{5ex}

\section{Relative presentations}\label{sec:relpres}

A {\em relative presentation} $\hat P=\langle H, {\bf x}\mid {\bf \hat r} \rangle$ consists of a group $H$, a generating set ${\bf x}=\{ x_1,...,x_n\}$ and a relator list ${\bf \hat r}=\{ \hat r_1,...,\hat r_m\}\subseteq H*F({\bf x})$. 
In this section we compare diagrammatic reducibility of relative presentations to directed diagrammatic reducibility. The concept of asphericity for relative presentations was introduced by Bogley and Pride \cite{BP92} and expressed in the language of pictures rather than diagrams (combinatorial maps). Spherical diagrams and spherical pictures are dual concepts. All definitions concerning pictures are presented in detail in Section 3.1 of the survey article Bogley-Edjvet-Williams \cite{BEW18}. 

\begin{rmk} {\em In the definition of pictures as given in \cite{BP92} and  \cite{BEW18} arcs are labeled by generators, but discs do not carry relator labels. This can lead to confusion, for example if ${\bf r}$ contains two words that are cyclic permutations of each other, and it does become a problem if $\hat r_i=\hat r_j$ in the list ${\bf r}$, something we do not want to rule out (topologically this means that a 2-cell might not be determined by its attaching map). For that reason we label discs in pictures by relators. If a disc carries the relator label $\pm r$ then reading off the arc labels around that disc (clockwise $+$, anticlockwise $-$) gives a cyclic permutation of $r$ (taking corners into consideration). The two discs involved in a dipole are required to carry relator labels $+r$ and $-r$. }
\end{rmk}

A connected spherical picture is assumed to contain a disc (and does not just consist of a floating circle). A connecting arc in a dipole is any arc that connects the two discs of the dipole. 

Given a relative presentation $\hat P=\langle H, {\bf x}\mid {\bf \hat r} \rangle$. Choose a presentation $P_0=\langle {\bf x_0}\ |\ {\bf r_0}\rangle$ of the group $H$.
If $\hat r=x_{1}h_1x_{2}h_2\ldots x_{t}h_t\in {\bf \hat r}$, where $x_{i}\in {\bf x}^{\pm 1}$ and $h_i\in H$, choose words $u_i$ in ${\bf x_0}^{\pm 1}$, such that $u_i$ represents $h_i$ in $G(P_0)=H$. Let $r=x_{1}u_1x_{2}u_2\ldots u_{t}$ and let ${\bf r}$ be the set of words built in this fashion from the set $\hat {\bf r}$. Let $P=\langle {\bf x_0}\cup {\bf x}\ | {\bf r_0}\cup {\bf r}\rangle$. A connected spherical picture $\mathbb P$ over $P$ contains a dipole with connecting arc labeled by $x$ if and only if the corresponding spherical diagram (obtained by dualizing) $f_{\mathbb P}\colon C\to K(P)$ contains a folding edge labeled by $x$. Thus we have the following reformulation of Definition \ref{defDrdira} in terms of pictures.

\begin{defn} Let $P=\langle {\bf x_0}\cup{\bf x}\ |\ {\bf r}\rangle$ be a presentation. Then $P$ is {\em DR directed away from ${\bf x_0}$} if every connected spherical picture that contains an arc with label from ${\bf x}$ also contains a dipole with a connecting arc with label from ${\bf x}$.
\end{defn}

\begin{lemma}\label{lemma:shrinkblow} Assume $H\to G(P)$ is injective.
\begin{enumerate}
\item  If there exists a connected spherical picture $\mathbb P$ over $P$ that contains $x\in {\bf x}$ as an arc label, and all connecting arcs of dipoles carry labels from ${\bf x_0}$, then there exists a reduced  connected spherical picture $\hat{\mathbb P}$ over $\hat P$ that contains $x$ as an arc label. 
\item Conversely, a reduced connected spherical picture $\hat{\mathbb P}$ over $\hat P$ that carries the label $x$ gives rise to a connected spherical picture 
$\mathbb P$ over $P$ that carries the label $x$ and all connecting arcs of dipoles (in case there are any) carry labels from ${\bf x_0}$.
\end{enumerate}
\end{lemma}

\noindent Proof. Suppose $\mathbb P$ is a connected spherical picture over $P$ as in the statement of (1). We remove all discs from $\mathbb P$ that come from relators in $\bf{r_0}$ and all arcs that carry labels from $\bf{x_0}$ to obtain a spherical picture $\mathbb P'$. Let $D_i$ be a disc in $\mathbb P$ associated with a relator $r=x_{1}u_1x_{2}u_2\ldots u_{t}\in{\bf r}$. The corresponding disc $D_i'\in\mathbb P'$ has a corner for every $u_i$, and we write the group element $h_i\in H$ represented by $u_i$ into that corner. We choose a connected component $\hat{\mathbb P}$ from $\mathbb P'$ that contains a disc with the arc label $x$ (such a disc exists because we assumed it existed in $\mathbb P$). We have to check that $\hat{\mathbb P}$ is indeed a spherical picture over the relative presentation $\hat P$. So suppose $\Delta$ is an inner region in $\hat{\mathbb P}$ and let $h_1...h_l$ be the sequence of corner labels encountered when reading along a simply closed loop $\gamma$ close to the boundary of $\Delta$. Reading along $\gamma$ in the original picture $\mathbb P$ gives a word $u=u_1...u_l$ in ${\bf{x_0}}^{\pm 1}$. This word $u$ represents the trivial element in $G(P)$ because it is the boundary word of a picture over $P$. Since we assumed that $G(P_0)\to G(P)$ is injective it follows that $u$ represents the trivial element in $G(P_0)=H$. So $h_1...h_l=1$. Finally note that $\hat{\mathbb P}$ is reduced. Otherwise it would contain a dipole which would result in a dipole in $\mathbb P$ with a connecting arc with label in $\bf x$, contradicting our assumption.

For the proof of (2) assume $\hat{\mathbb P}$ is a spherical picture over $\hat P$ as in the statement. If $D_i$ is a disc in $\hat{\mathbb P}$ associated $\hat r=x_{1}h_1x_{2}h_2\ldots x_{t}h_t\in {\bf \hat r}$, then add a sequence of arc segments to the corner labeled $h_i$, so that reading across the sequence gives the word $u_i$. From every inner region $\Delta$ remove a small open disc $d$ from the interior. Lead all $u$ arc sequences coming from the corners of $\Delta$ to the boundary of $d$. Note that the word $w_d$ we read off the boundary of $d$ represents the trivial element of $G(P_0)=H$ since $\hat{\bf{\mathbb P}}$ is a picture over $\hat P$. Thus there exists a picture $\mathbb P_d$ over $P_0$ with boundary word $w_d$. We insert that picture as a replacement for $d$. For the outer region lead all $u$ arc sequences to the boundary of the ambient disc $D$ of the original picture $\hat{\mathbb P}$. What we have built is a picture $\mathbb P_1$ over $P$ whose boundary is a word $w$ in ${\bf{x_0}}^{\pm 1}$. So $w$ represents the trivial element of $G(P)$. Since we assumed that $H\to G(\hat P)$ is injective it follows that $H=G(P_0)\to G(P)$ is injective. So $w$ represents the trivial element of $G(P_0)$ and it follows that there exists a picture $\mathbb P_w$ over $P_0$ with boundary $w$. We add $\mathbb P_w$ to $\mathbb P_1$ and obtain a spherical picture $\mathbb P_2$ over $P$. Let $\mathbb P$ be a 
connected component in $\mathbb P_2$ that contains a disc with arc label $x$. It exists because we assumed it exists in $\hat{\mathbb P}$. Note that $\mathbb P$ can contain a dipole. But since all the arcs we added have labels from $\bf{x_0}$, all connecting arcs of dipoles have labels from $\bf{x_0}$. \qed \\

Again we start with a relative presentation $\hat P=\langle H, {\bf x}\ |\ {\bf \hat r} \rangle$ and built $P=\langle {\bf x\cup x_0 }\ |\ {\bf r_0}\cup {\bf r} \rangle$ as above. A relative presentation $\hat P$ is called DR if every connected spherical picture over $\hat P$ contains a dipole (see Definition 3.1 in \cite{BEW18}).

\begin{thm}\label{thm:equivalence} $\hat P$ is DR if and only if $P$ is DR directed away from ${\bf x_0}$.
\end{thm}

\noindent Proof. Assume that the relative presentation $\hat P$ is DR. It follows that $H\to G(\hat P)$ is injective (see Theorem 3.4 in \cite{BEW18}). Let $\mathbb P$ be a spherical picture that carries an edge label $x\in {\bf x}$, and all connecting arcs of dipoles (should there be any) carry labels from ${\bf x_0}$. By Lemma \ref{lemma:shrinkblow} (1) there exists a reduced connected spherical picture $\hat{\mathbb P}$ that carries the edge label $x$. This contradicts the assumption that $\hat P$ is DR. 

For the other direction, assume that $P$ is DR directed away from ${\bf x_0}$. Assume we have a reduced spherical diagram $\hat{\mathbb P}$ over $\hat P$. Since it is not empty it must carry a label $x\in {\bf x}$. Since $P$ is DR directed away from ${\bf x_0}$ we know from Theorem \ref{RelPiS} (2) that $G(P_0)=H$ injects into $G(P)$ and thus $H$ injects into $G(\hat P)$. By Lemma \ref{lemma:shrinkblow} (2) there exists a spherical picture that carries the label $x\in {\bf x}$ whose dipoles (should there be such) all have connecting arcs that carry labels from ${\bf x_0}$. This contradicts the assumption that $P$ is DR away from ${\bf x_0}$. \qed \\

We note that Theorem \ref{thm:equivalence} implies the following: If $P_1=\langle {\bf x\cup x_1 }\ |\ {\bf r_1}\cup {\bf r} \rangle$ and $P_2=\langle {\bf x\cup x_2 }\ |\ {\bf r_2}\cup {\bf r} \rangle$ are two presentations built from the relative presentation $\hat P=\langle H, {\bf x}\ |\ {\bf \hat r} \rangle$, where $\langle {\bf x_1 }\ |\ {\bf r_1} \rangle$ and $\langle {\bf x_2 }\ |\ {\bf r_2} \rangle$ are presentations of $H$, then $P_1$ is DR directed away from ${\bf x_1}$ if and only if $P_2$ is DR directed away from ${\bf x_2}$.\\

Any result about directed DR has a translation into a result about DR for relative presentations and vice versa. We will next give a version of the generalized Corson-Trace Theorem \ref{thm:GenCT} for relative presentations. We start with a relative presentation $\hat P=\langle H, {\bf x}\ |\ {\bf \hat r} \rangle$ and built $P=\langle {\bf x\cup x_0 }\ |\ {\bf r_0}\cup {\bf r} \rangle$. 
Let $p\colon \tilde K(P) \to K(P)$ be the universal covering. We contract each component of $p^{-1}(K(P_0))$ to a point and arrive at a simply-connected 2-complex $Y$. Denote the quotient map by $q\colon \tilde K(P) \to Y$.

\begin{thm}\label{thm:CTrelP} The relative presentation $\hat P$ is DR if and only if every finite subcomplex of $Y$ collapses into $Y^{(1)}$.
\end{thm} 

\begin{rmk}\em We have not found this result explicitly  stated in the literature but Theorem 3.6 in Bogley-Edjvet-Williams \cite{BEW18} is closely related. The proof given there does use the complex $Y$ (called $X$ there) and the quotient map $q\colon \tilde K(P) \to Y$ and it is argued that the assumption $\hat P$ is DR implies that $Y$ is DR. So by the classical Corson-Trace every finite subcomplex of $Y$ collapses into the 1-skeleton. We note that the implication from (b) to (c) in the given proof does seem to require that the relative presentation presents a finite group.
\end{rmk}

\noindent Proof of Theorem \ref{thm:CTrelP}. By Theorem \ref{thm:GenCT} and Theorem \ref{thm:equivalence} we know that $\hat P$ is DR if and only if every finite subcomplex of $\tilde K(P)$ collapses into $p^{-1}(K(P_0))\cup \tilde K(P)^{(1)}$ where each collapse uses a free edge of the form $(g,x)$, $x\in {\bf x}, g\in G(P)$. Note also that the quotient map $q\colon \tilde K(P) \to Y$ maps 2-cells that are not contained in $p^{-1}(K(P_0))$ bijectively to the 2-cells in $Y$. In particular, if $(g,x)$, $x\in {\bf x}, g\in G(P)$ is a free edge of a subcomplex $M$ of $\tilde K(P)$, then $q(g,x)$ is a free edge in the subcomplex $q(M)$.

Suppose first that $\hat P$ is DR. Let $M'$ be a finite subcomplex of $Y$. Then there exists a finite subcomplex $M$ of $\tilde K(P)$ such that $q(M)=M'$ ($q^{-1}(M')$ might not be finite, but it is finite outside of $p^{-1}(K(P_0))$). Now $M$ collapses into $p^{-1}(K(P_0))\cup \tilde K(P)^{(1)}$, where each collapse uses a free edge of the form $(g,x)$, $x\in {\bf x}, g\in G(P)$. By as we just mentioned, each $q(g,x)$ is a free edge in $q(M)=M'$. Thus the collapsing sequence can be applied to $q(M)=M'$, and it collapses $M'$ into $Y^{(1)}$.

For the other direction assume that $M$ is a finite subcomplex of $\tilde K(P)$. Then $q(M)=M'$ is a finite subcomplex of $Y$ and it follows that $M'$ collapses into $Y^{(1)}$. Let $(g,x)$, $x\in {\bf x}, g\in G(P)$, be an edge in $M$ that maps to a free edge in $q(M)=M'$. But then $(g,x)$ is a free edge in $M$, by the remark in the first paragraph of the proof. \qed \\

We end this section with a simple lemma that will be useful for establishing directed DR in concrete situations. Let $\phi\colon H_1\to H_2$ be a group homomorphism and suppose $\hat P_1=\langle H_1, {\bf x}\mid {\bf \hat r}\rangle $ is a relative presentation. If $\hat r=x_1h_1x_2h_2\ldots x_kh_k$, where the $x_i\in {\bf x}^{\pm 1}$, $h_i\in H$, define $\phi(\hat r)=x_1\phi(h_1)x_2\phi(h_2)\ldots x_k\phi(h_k)$. It could happen that there are relators $\hat r_1$ and $\hat r_2$ in $\hat {\bf r}$ such that $\phi(\hat r_1)=\phi(\hat r_2)$. In that case we would consider $\phi(\hat r_1)$ and $\phi(\hat r_2)$ distinct element on the list $\phi({\bf r})$.  So $\phi\colon \hat{\bf r}\to \phi(\hat{\bf r})$ is a bijection of lists. We have a relative presentation $\hat P_2=\langle H_2, {\bf x}\mid \phi({\bf \hat r})\rangle$ and write 
$$\phi\colon \hat P_1=\langle H_1, {\bf x}\mid {\bf \hat r}\rangle  \to \hat P_2=\langle H_2, {\bf x}\mid \phi({\bf \hat r})\rangle.$$

\begin{lemma}\label{lemma:phi} Suppose we have $\phi\colon \hat P_1 \to \hat P_2$ induced by a group homomorphism $H_1\to H_2$.  If $\hat P_2$ is DR, then $\hat P_1$ is so as well.
\end{lemma}

\noindent Proof. Note that $\phi$ takes spherical pictures over $\hat P_1$ to spherical pictures over $\hat P_2$. In fact if $\hat{\mathbb P}$ is a connected spherical picture over $\hat P_1$, then $\phi(\hat{\mathbb P})$ looks the same, except that relator labels change from $\hat r$ to $\hat \phi(r)$ and corner label change from $h$ to $\phi(h)$ in $\phi(\hat{\mathbb P})$. 

Let $\hat{\mathbb P}$ be a connected spherical picture over $\hat P_1$. Then $\phi(\hat{\mathbb P})$ is a connected spherical picture over $\hat P_2$ and hence contains a dipole with label $\pm \phi(\hat r)$, say. Then $\hat{\mathbb P}$ contains a corresponding dipole with label $\pm\hat r$. \qed

\vspace{5ex}

\pagebreak

\section{Tests for directed DR}\label{sec:tests}

Let $T_1=\langle x_1,\ldots ,x_k\mid r_1,\ldots ,r_l\rangle$ and $T_2=\langle y_1,\ldots ,y_p\mid s_1,\ldots ,s_q\rangle$ be presentations. Suppose we have a map $\phi_0\colon F( x_1,\ldots ,x_k)\to F(y_1,\ldots ,y_p)$ that induces a group homomorphism $G(T_1)\to G(T_2)$. For $n>k$ we can extend $\phi_0$ to $\phi\colon F(x_1,\ldots ,x_n)\to F(y_1,\ldots ,y_p, x_{k+1},\ldots ,x_n)$ by defining $\phi(x_i)=x_i$ for $i>k$. Now let $P_1$ be a presentation of the form 
$$P_1=\langle x_1,\ldots ,x_k, x_{k+1},\ldots ,x_n\mid r_1,\ldots ,r_l, r_{l+1},\ldots ,r_m\rangle.$$ We assume each relator $r_j$, $j>l$, contains a generator $x_i$, $i>k$, so that $P_{S_1}=T_1$ for $S_1=\{ x_1,\ldots ,x_k\}$. $$P_2=\langle y_1,\ldots ,y_p, x_{k+1},\ldots ,x_n \mid s_1,\ldots ,s_q, \phi(r_{l+1}),\ldots ,\phi(r_m)\rangle.$$ We use $\phi\colon P_1\to P_2$ as shorthand for the situation just described. Here is Lemma \ref{lemma:phi} in the directed DR setting.

\begin{lemma}\label{sbarP} Suppose we have $\phi\colon P_1 \to P_2$. If $P_2$ is DR directed away from $S_2=\{ y_1,\ldots ,y_p\}$, then $P_1$ is DR directed away from $S_1=\{ x_1,\ldots ,x_k\}$.
\end{lemma}

\noindent Proof. Let 
$$\hat P_1=\langle G(T_1), x_{k+1},...,x_n\ |\ \hat{r}_{l+1},\ldots ,\hat r_m\rangle$$ and 
$$\hat P_2=\langle G(T_2), x_{k+1},\ldots ,x_n \mid \hat \phi(r_{l+1}),\ldots ,\hat \phi(r_m)\rangle$$ be the associated relative presentations. Note that we have $\phi_0\colon \hat P_1\to \hat P_2$ as in Lemma \ref{lemma:phi}. 
Since $P_2$ is DR directed away from $S_2$ it follows from Theorem \ref{thm:equivalence} that $\hat P_2$ is DR. It follows from Lemma \ref{lemma:phi} that $\hat P_1$ is DR. Hence, by  Theorem \ref{thm:equivalence} $P_1$ is DR directed away from $S_1$.\qed \\

The simplest choice for $T_2$ is the empty presentation $T_2=\langle \ |\ \rangle$, that is $G(T_2)$ is trivial as in the following example.

\begin{example} \em Consider $$P_1=\langle x_1,\ldots ,x_k, a, b\mid u_1au_2bu_3a^{-1}u_4b^{-1} \rangle \to P_2= \langle a, b\mid aba^{-1}b^{-1}\rangle,$$
where the $u_i$ are words in $x_1^{\pm 1},\ldots ,x_k^{\pm 1}$. Since $P_2$ is $DR$ (directed away from $\emptyset$), it follows that $P_1$ is DR directed away from $S=\{ x_1,\ldots ,x_k \}$. More general: Take any DR presentation $P_2$. Add generators $S=\{ x_1,\ldots ,x_k \}$ and insert words in $S^{\pm 1}$ into the relators of $P_2$ and one obtains a presentation $P_1$ which is DR directed away from $S$.
\end{example} 

A transformation of spherical diagrams over $K(P_1)$ to diagrams over $K(P_2)$ can be seen directly. For every relation $r$ in $T_1$ choose a Van Kampen diagram $M_{r}\to K(T_2)$ that expresses $\phi_0(r)$ in terms  of the $s$-relations in $T_2$. This is possible because we assumed that $\phi_0$ induces a group homomorphism $G(T_1)\to G(T_2)$. Now given a spherical diagram $f\colon C\to K(P_1)$ we can produce a diagram $\phi\circ f: C'\to K(P_2)$ in the following way. If $e$ is an edge in $C$ labeled with $x$, then subdivide it and label it with $\phi(x)$. Note that if $x$ is a generator of $T_1$, then $\phi(x)$ is a word in the $y$ generators of $T_2$. If $d$ is a 2-cell in $C$ labeled with an $r_i$,
$1\le i\le k$  from $T_1$, then remove the interior of $d$ and insert $M_{r_i}$. If $d$ is a 2-cell in $C$ labeled with $r_i$, $k+1\le i \le n$, then insert a relator disc labeled $\phi(r_i)$. The resulting diagram over $K(P_2)$ is $\phi\circ f$. Using this transformation we can give a direct proof of Lemma \ref{sbarP} without reference to relative presentations. If $f\colon C\to K(P_1)$ is a spherical diagram that contains an edge labeled $x_i$, $k+1\le i\le n$ then so does $\phi\circ f$. Since $P_2$ is DR away from $S_2$, $\phi\circ f$ contains a folding edge with label $x_j$, $k+1\le j\le n$ (not in $S_2$). This edge is present in $f$ and is a folding edge.\\

A {\it cycle} in a graph $\Gamma $ is a closed edge path in $\Gamma $. If $P=\langle x_1,\ldots ,x_n\mid r_1,\ldots ,r_m\rangle$ is a presentation then the {\it Whitehead graph} $W(P)$ is the boundary of a regular neighborhood of the only vertex of $K(P)$.  It is a non-oriented graph on vertices $\{ x_i^+, x_i^- \mid i=1,\ldots ,n\}$, where $x_i^+$ is a point of the oriented edge $x$ of $K(P)$ close to the beginning of that edge, and $x_i^-$ is a point close to the ending of that edge. The {\em positive graph} $W^+(P)\subset W(P)$  is the full subgraph on the vertex set  $\{ x_i^+\mid i=1,\ldots ,n\}$, the {\em 
negative graph } $W^-(P)\subset W(P)$ is the full subgraph on the
vertex set $\{x_i^- \mid i=1,\ldots ,n\}$. Let $C$ be a cell decomposition of a 2-sphere with oriented edges. A {\em sink} is a vertex in $C$ with all adjacent edges pointing towards it, a {\em source} is a vertex in $C$ with all adjacent edges pointing away from it.

\begin{thm}\label{sminmax} Let $Q=\langle x_{1},\ldots ,x_n\mid r_{1},\ldots ,r_m\rangle$ be a finite presentation with cyclically reduced relators of exponent sum 0. Assume that $W^+(Q)$ or $W^-(Q)$ is a forest. Then $Q$ is DR directed away from $\{ x_k\}$ for each generator $x_k$.
\end{thm}

\noindent Proof. Assume that $W^-(Q)$ is a forest. Let $f\colon C\to K(Q)$ be a spherical diagram. Since the relators are assumed to be of exponent sum zero, it follows that $C$ has a sink and a source vertex. Suppose $v$ is a sink vertex. The link at $v$ gives a cycle $z$  in a connected component of $W^-(Q)$, which is a tree. Thus, the edges covered by $z$ form a subtree $\Gamma$ of $W^-(Q)$, which contains two (or more) distinct vertices of valency one of the set $\{ x_{1}^-,\ldots ,x_n^-\}$. Thus the edges at $v$ in $C$, all directed towards $v$, will contain folding edges labeled by two elements of $\{ x_{1},\ldots ,x_n\}$. Thus one of them is labeled by a generator distinct from $x_k$. This implies that $Q$ is DR directed away from $\{ x_k\}$.

Analogous arguments apply if $W^+(Q)$ is a forest for a source vertex of $C$. \qed \\

It was known previously (Gersten  \cite{Ger87}) that presentations satisfying the assumptions of Theorem \ref{sminmax} are DR. Lemma \ref{sbarP} together with Theorem \ref{sminmax} provide a tool for showing directed DR for presentations. We will discuss this in greater detail in the next section.\\

If $w$ is a real valued function on the edges of a graph $\Gamma $ and $z=e_1\ldots e_p$ is an edge path in $\Gamma$, we write 
$$w(z)=\sum_{1\le i\le p}w(e_i).$$

The following is a weight test for directed diagrammatic reducibility. It is a generalized version of Gerstens weight test (see \cite{Ger87}).

\begin{thm}\label{sWT} Let $Q=\langle y_1,\ldots ,y_p,x_1,\ldots ,x_n\mid r_1,\ldots ,r_m\rangle$ be a presentation with cyclically reduced relators.
Suppose we can assign non-negative weights $\omega(e)\ge 0$ to the edges $e$ in $W=W(Q)$, such that the following conditions hold: 
\begin{enumerate}
\item If $e$ is an edge of $W$ that connects (for some $i,j\le p$) $y_i^\epsilon$ with $y_j^\delta$, ($\epsilon ,\delta =\pm$) then $\omega(e)\ge 1$; 
\item If $e$ is an edge of $W$ and exactly one of its boundary points is $y_i^+$ or $y_i^-$ (for some $i\le p$), then $\omega(e)\ge 1/2$;
\item If $z$ is a reduced cycle in $W$, then $\omega(z)\ge 2$;
\item Let $r$ be a relator of length $d(r)$ from ${Q}$, then $\sum_{c\in r}\omega(c)\le d(r)-2$.
\end{enumerate}
Then $Q$ is DR directed away from $\{ y_1,\ldots ,y_p\}$.
\end{thm}

\noindent Proof. Note that conditions (3) and (4) in the statement of the theorem imply that $Q$ satisfies the standard weight test and it follows that $Q$ is DR. Let $h\colon C\to K(Q)$ be a spherical diagram. DR implies that there are folding edges in $C$. Suppose only the generators $y_1,\ldots ,y_p$ occur as labels on folding edges.

Let $z$ be a cycle in $W$ obtained from reading around a vertex in $C$. If $z$ is reduced then $\omega(z)\ge 2$ by condition (3). Suppose $z$ is not reduced. Then $z$ contains consecutive reducing edge pairs $e^{\epsilon}e^{-\epsilon}$, where $\epsilon=\pm 1$. For convenience we assume $\epsilon=1$.  The terminal vertex of $e$ must be $y_i^+$ or $y_i^-$ for some $y_i$. If the other vertex of $e$ (note that $e$ connects two distinct vertices in $W$ because we assume that relators in $P$ are cyclically reduced) comes also from a generator $y_j$, then the pair $ee^{-1}$ contributes at least $2$ to the weight of $z$, by condition (1). Since weights are positive we then have $\omega(z)\ge 2$. 

So suppose the other vertex of $e$ does come from a generator $x_i$. Then, by condition (2), the edge pair $ee^{-1}$ contributes at least $1$ to the weight of $z$. So if $z$ contains a second such edge pair, then the two edge pairs contribute at least $2$ to the weight of $z$, and we again have $\omega(z)\ge 2$. Suppose that $z$ contains only one consecutive reducing edge pair $ee^{-1}$.
Then there is a reduced cycle $z'$ which consists of edges which also appear in $z$. By condition (3) we have $\omega(z)\ge \omega(z')\ge 2$.

Since $\omega(z)\ge 2$ for every cycle in $W$ obtained by reading around a vertex in $C$, the combinatorial curvature at every vertex in the spherical diagram $C$ is less or equal to $0$. The combinatorial curvature of the 2-cells in $C$ is also less or equal to $0$ by condition (4). The combinatorial Gauss-Bonnet Theorem now implies that the Euler characteristic of $C$ is less or equal to $0$, which contradicts the fact that $C$ is a sphere. So in $h$ must be a folding edge labeled by a generator different from one of $\{ y_1,\ldots ,y_p\}$.\qed \\

Again, Lemma \ref{sbarP} together with Theorem \ref{sWT} provide a method for proving directed DR for presentations. Details are discussed in the next section.

\vspace{5ex}

\section{Applications}\label{sec:appl}

The small cancellation conditions C($p$), T($q$) are defined for instance in the book of Lyndon and Schupp (see \cite{LS77}). The following is a generalization of the observation of Gersten that small cancellation groups are DR (see Remark 4.18 of \cite{Ger87}).

\begin{thm}\label{s44} Let $Q$ be a finite presentation with cyclically reduced relators and $S$ a subset of the generators. Assume that $Q$ is C(4), T(4) or C(6), T(3) and that no two consecutive letters in a (cyclically read) relator of $Q$ are elements of $S^{\pm 1}$. Then $Q$ is DR directed away from $S$.
\end{thm}

\noindent Proof. We define a weight function $w\colon E\to\mathbb R$ where $E$ are the edges of the Whitehead graph $W(Q)$, such that the weight test of Theorem \ref{sWT} is satisfied with $S=\{ y_1,\ldots ,y_p\}$.

Give each edge $e\in E$, which is in a cycle of length 2 weight $1$. Give all other edges in case C(4), T(4) weight $1/2$ and in case C(6), T(3) weight $2/3$. So condition (2) of Theorem \ref{sWT} is satisfied. T(4) (or T(3)) implies condition (3) and C(4) (or C(6)) implies condition (4).
Since no two consecutive letters in a (cyclically read) relator of $Q$ are generators of $S$ there will be no edges in $E$ connecting $y_i^\pm $ with $y_k^\pm $ for some $y_i,y_k\in S$. So condition (1) is satisfied also.\qed \\

At last we present some applications to labeled oriented trees. Labeled oriented trees were introduced by J. Howie in connection with the Whitehead Conjecture (see \cite{How85}).
A {\em labeled oriented tree} (LOT) $\mathcal P$ consists of an oriented graph which is a tree, whose oriented edges are labeled by its vertices. From a labeled oriented tree  $\mathcal P$ we obtain a LOT presentation $P$: The generators in $P$ are the vertices of $\mathcal P$, and for every oriented edge from a vertex $u_1$ to a vertex $u_2$, labeled by $u_3$ we have a relation $u_1u_3(u_3u_2)^{-1}$. Wirtinger presentations read from knot diagrams are examples of LOT presentations. Much information on labeled oriented trees and their importance to questions in combinatorial topology and group theory can be found in Rosebrock \cite{Ro18}. In \cite{HaRo17} the authors show the combinatorial asphericity of injective LOTs (i.e. LOTs where different edges carry different edge labels).

In Rosebrock \cite{Ro94} it is described how to check whether a LOT is C(4), T(4).
If a LOT is C(4), T(4) then Theorem \ref{s44} implies that it is DR away from any of its generators. But there is more:

\begin{example} Consider the LOT $Q$ of Figure \ref{alc4t4} with any orientation of its edges. This LOT is C(4), T(4). If you choose $S$ to be one of the sets
\[\{ x_1,x_2,x_5\} , \{ x_1,x_2\} , \{ x_1,x_4\} , \{ x_1,x_5\} , \{ x_2,x_5\} , \{ x_2,x_6\} , \{ x_2,x_7\} , \{ x_3,x_6\} , \{ x_3,x_7\} , \{ x_4,x_7\} \]
then Theorem \ref{s44} implies that $Q$ is DR away from $S$.
\end{example}
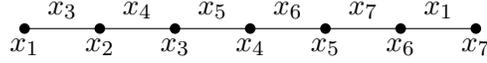
\begin{figure}[ht]\centering
\begin{tikzpicture}
\draw (1,0)--(2,0)--(3,0)--(4,0)--(5,0)--(6,0)--(7,0);
\fill (1,0) circle (2pt);
\fill (2,0) circle (2pt);
\fill (3,0) circle (2pt);
\fill (4,0) circle (2pt);
\fill (5,0) circle (2pt);
\fill (6,0) circle (2pt);
\fill (7,0) circle (2pt);
\node [below] at (1,0) {$x_1$};
\node [above] at (1.5,0) {$x_3$};
\node [below] at (2,0) {$x_2$};
\node [above] at (2.5,0) {$x_4$};
\node [below] at (3,0) {$x_3$};
\node [above] at (3.5,0) {$x_5$};
\node [below] at (4,0) {$x_4$};
\node [above] at (4.5,0) {$x_6$};
\node [below] at (5,0) {$x_5$};
\node [above] at (5.5,0) {$x_7$};
\node [below] at (6,0) {$x_6$};
\node [above] at (6.5,0) {$x_1$};
\node [below] at (7,0) {$x_7$};
\end{tikzpicture}
\caption{\label{alc4t4} A labeled oriented tree which is C(4), T(4).}
\end{figure}
 
A {\em sub-LOT} of a labeled oriented tree $\mathcal P$ is a connected subgraph $\mathcal T$ (containing at least one edge) such that each edge label of $\mathcal T$ is a vertex of $\mathcal T$. A LOT is called {\it compressed} if each relator consists of three different generators. The relators of a LOT presentation of a compressed LOT are cyclically reduced.

 Let $\mathcal P$ be a compressed LOT and $\mathcal T$ be a maximal proper sub-LOT with vertex set $S$. Let $\bar{\mathcal P}$ be the LOT obtained from $\mathcal P$ by collapsing all of $\mathcal T$ to a vertex $y$. Every occurrence of a vertex $x$ of $\mathcal T$ in $\mathcal P-\mathcal T$ is replaced by $y$ in $\bar{\mathcal P}$. Let $T$, $P$, $\bar P$ be the corresponding LOT presentations.

\begin{thm}\label{sLOT} If $\bar P$ is compressed and DR directed away from $y$ then $P$ is DR directed away from $S$, the generators of $T$. In particular if
\begin{enumerate}
\item $W^+(\bar P)$, or $W^-(\bar P)$ is a tree, or
\item the Whitehead graph $W(\bar P)$ does not contain cycles of length less than four,
\end{enumerate}
then $P$ is DR directed away from $S$.
\end{thm}

\noindent Proof. Let $\bar T=\langle y \mid - \rangle$. We can apply Lemma \ref{sbarP} to the situation $T\subseteq P \to \bar T\subseteq \bar P$, where all generators of $T$ go to $y$ in $\bar T$. Here are the details. Let $S=\{ x_1,\ldots , x_k\}$ be the generators of $T$ (the vertices in $\mathcal T$) and let $\{ x_1,\ldots ,x_n\}$ be the generators of $P$ (the vertices in $\mathcal P$). Define $\phi_0\colon F(x_1,\ldots ,x_k)\to F(y)$ to be the map that sends each $x_i$ to $y$, and $\phi\colon F(x_1,\ldots ,x_n)\to F(y, x_{k+1},\ldots ,x_n)$ the extension of $\phi_0$ that sends $x_i$ to $x_i$ for $k+1\le i\le n$. Then 
$$\langle y, x_{k+1},\ldots ,x_n \mid \phi(r_{k}),\ldots , \phi(r_{n-1})\rangle$$ 
is the LOT presentation $\bar P$. Here the $r_j$, $k\le j \le n-1$ are the LOT relations coming from the edges in $\mathcal P-\mathcal T$.

Now (1) implies that $\bar P$ is DR away from $\{ y\}$ by Theorem \ref{sminmax} and (2) implies that $\bar P$ is DR away from $\{ y\}$ by Theorem \ref{s44}. It follows from Lemma \ref{sbarP} that $P$ is DR directed away from $S$. \qed \\

\begin{example} {\em Figure \ref{aex1} below shows a compressed LOT $\mathcal P$ with a sub-LOT $\mathcal T$ with vertices $S=\{ x_1,\ldots ,x_5\}$. Below $\mathcal P$ we see the LOT $\bar{\mathcal P}$ obtained from $\mathcal P$ by collapsing $\mathcal T$ to the vertex $y$.}\\
\end{example}

\vspace{-4ex}

\begin{figure}[ht]\centering
\begin{tikzpicture}
\node at (-1,0) {$\mathcal P$};
\node[red] at (5,0.75) {$\mathcal T$};

\fill (0,0) circle (2pt);
\node[below] at (0,0) {$u_3$};
\fill (1,0) circle (2pt);
\node[below] at (1,0) {$u_2$};
\fill (2,0) circle (2pt);
\node[below] at (2,0) {$u_1$};
\fill[red] (3,0) circle (2pt);
\node[red,below] at (3,0) {$x_1$};
\fill[red] (4,0) circle (2pt);
\node[red,below] at (4,0) {$x_2$};
\fill[red] (5,0) circle (2pt);
\node[red,below] at (5,0) {$x_3$};
\fill[red] (6,0) circle (2pt);
\node[red,below] at (6,0) {$x_4$};
\fill[red] (7,0) circle (2pt);
\node[red,below] at (7,0) {$x_5$};
\fill (8,0) circle (2pt);
\node[below] at (8,0) {$u_4$};

\begin{scope}[decoration={markings, mark=at position 0.5 with {\arrow{>}}}]
\draw [red,postaction={decorate}] (3,0) -- (4,0) node[midway, above]{$x_3$};
\draw [red,postaction={decorate}] (4,0) -- (5,0) node[midway, above]{$x_1$};
\draw [red,postaction={decorate}] (5,0) -- (6,0) node[midway, above]{$x_5$};
\draw [postaction={decorate}] (7,0) -- (8,0) node[midway, above]{$u_3$};
\end{scope}

\begin{scope}[decoration={markings, mark=at position 0.5 with {\arrow{<}}}]
\draw [postaction={decorate}] (0,0) -- (1,0) node[midway, above]{$u_1$};
\draw [postaction={decorate}] (1,0) -- (2,0) node[midway, above]{$u_3$};
\draw [postaction={decorate}] (2,0) -- (3,0) node[midway, above]{$u_4$};
\draw [red,postaction={decorate}] (6,0) -- (7,0) node[midway, above]{$x_1$};
\end{scope}

\begin{scope}[yshift=-1cm]
\node at (-1,0) {$\bar{\mathcal P}$};

\fill (0,0) circle (2pt);
\node[below] at (0,0) {$u_3$};
\fill (1,0) circle (2pt);
\node[below] at (1,0) {$u_2$};
\fill (2,0) circle (2pt);
\node[below] at (2,0) {$u_1$};
\fill[red] (3,0) circle (2pt);
\node[red,below] at (3,0) {$y$};
\fill (4,0) circle (2pt);
\node[below] at (4,0) {$u_4$};

\begin{scope}[decoration={markings, mark=at position 0.5 with {\arrow{>}}}]
\draw [postaction={decorate}] (3,0) -- (4,0) node[midway, above]{$u_3$};
\end{scope}

\begin{scope}[decoration={markings, mark=at position 0.5 with {\arrow{<}}}]
\draw [postaction={decorate}] (0,0) -- (1,0) node[midway, above]{$u_1$};
\draw [postaction={decorate}] (1,0) -- (2,0) node[midway, above]{$u_3$};
\draw [postaction={decorate}] (2,0) -- (3,0) node[midway, above]{$u_4$};
\end{scope}

\end{scope}

\end{tikzpicture}
\caption{\label{aex1} The LOT $\mathcal P$ with sub-LOT $\mathcal T$ and the LOT $\bar{\mathcal P}$}
\end{figure}
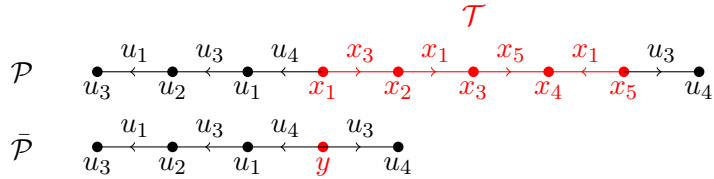

\noindent Note that $W^+(\bar{P})$ is a tree, so by Theorem \ref{sLOT} $P$ is DR directed away from $S$. Also observe that $W^-(T)$ is a tree, so $K(T)$ is DR, which implies that $K(P)$ is aspherical by Theorem \ref{RelPiS} (1). Note that neither $W^+(P)$ nor $W^-(P)$ is a tree.\\

The process of collapsing a sub-LOT in a given labeled oriented tree can also be reversed: If $\bar {\mathcal P}$ is a labeled oriented tree, $y$ is a vertex in $\bar {\mathcal P}$, and $\mathcal T$ is a labeled oriented tree, we remove $y$ from $\bar {\mathcal P}$ and insert $\mathcal T$ to obtain a labeled oriented tree $\mathcal P$ that contains $\mathcal T$. Collapsing $\mathcal T$ in $\mathcal P$ to a vertex $y$  brings us back to $\bar {\mathcal P}$. So the previous theorem can also be stated as follows: If $\bar{\mathcal P}$ is a labeled oriented tree that satisfies either condition (1) or (2) of Theorem \ref{sLOT}, then inserting any LOT $\mathcal T$ into $\bar{\mathcal P}$ results in a labeled oriented tree $\mathcal P$ for which the LOT presentation $P$ is DR directed away from the set $S$ of vertices of $\mathcal T$.\\

LOT presentations $P$ where $W^+({P})$ or $W^-({P})$ is a tree abound. 
If  a LOT $\mathcal P'$ is obtained from a LOT $\mathcal P$ by changing some edge orientations, we call $\mathcal P'$ a {\em reorientation} of $\mathcal P$. In \cite{HR01} Proposition 5.1, Huck and Rosebrock show that each LOT $\mathcal P$ has a reorientation $\mathcal P'$ such that its positive Whitehead graph $W^+(P')$ is a tree. Theorem \ref{sminmax} now implies

\begin{thm}
Each LOT $\mathcal P$ has a reorientation $\mathcal P'$ so that $P'$ is DR away from any one of its generators.
\end{thm}

In \cite{Ro94} Rosebrock gives conditions on a labeled oriented tree so that condition (2) of Theorem \ref{sLOT} holds. A concrete example is shown in Figure \ref{aex2}. 

\begin{example} {\em Figure \ref{aex2} below shows a labeled oriented tree $\mathcal P$ (orientations can be chosen at will) with a sub-LOT $\mathcal T$ between $u_4$ and $u_4'$ (which can be filled in at will).}
\end{example}

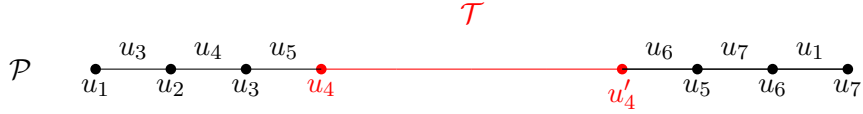
\begin{figure}[ht]\centering
\begin{tikzpicture}
\node at (-1,0) {$\mathcal P$};
\node[red] at (5,0.75) {$\mathcal T$};

\fill (0,0) circle (2pt);
\node[below] at (0,0) {$u_1$};
\fill (1,0) circle (2pt);
\node[below] at (1,0) {$u_2$};
\fill (2,0) circle (2pt);
\node[below] at (2,0) {$u_3$};
\fill[red] (3,0) circle (2pt);
\node[red,below] at (3,0) {$u_4$};
\fill(7,0)[red] circle (2pt);
\node[red, below] at (7,0) {$u'_4$};
\fill (8,0) circle (2pt);
\node[below] at (8,0) {$u_5$};
\fill (9,0) circle (2pt);
\node[below] at (9,0) {$u_6$};
\fill (10,0) circle (2pt);
\node[below] at (10,0) {$u_7$};

\draw [red,postaction={decorate}] (3,0) -- (4,0); 
\draw [red,postaction={decorate}] (4,0) -- (5,0);
\draw [red,postaction={decorate}] (5,0) -- (6,0);
\draw [postaction={decorate}] (7,0) -- (8,0);
\draw [postaction={decorate}] (8,0) -- (9,0);
\draw [postaction={decorate}] (9,0) -- (10,0);

\draw [postaction={decorate}] (0,0) -- (1,0) node[midway, above]{$u_3$};
\draw [postaction={decorate}] (1,0) -- (2,0) node[midway, above]{$u_4$};
\draw [postaction={decorate}] (2,0) -- (3,0) node[midway, above]{$u_5$};
\draw [red,postaction={decorate}] (6,0) -- (7,0);
\draw [postaction={decorate}] (7,0) -- (8,0) node[midway, above]{$u_6$};
\draw [postaction={decorate}] (8,0) -- (9,0) node[midway, above]{$u_7$};
\draw [postaction={decorate}] (9,0) -- (10,0) node[midway, above]{$u_1$};

\end{tikzpicture}
\caption{\label{aex2} A labeled oriented tree $\mathcal P$ with sub-LOT $\mathcal T$.}
\end{figure}

\noindent Notice that if we collapse the red sub-LOT $\mathcal T$ to the vertex $y=u_4$, we obtain a labeled oriented tree $\bar{\mathcal P}$ for which $W(\bar{\mathcal P})$ does not contain cycles of length less than four. It follows from Theorem \ref{sLOT} that $P$ is DR away from the vertex set of $\mathcal T$.



\vspace{3ex}

\noindent Jens Harlander\\
Department of Mathematics\\
Boise State University\\
Boise, ID 83725-1555\\
USA\\
email: jensharlander@boisestate.edu\\

\noindent Stephan Rosebrock\\
P{\"a}dagogische Hochschule Karlsruhe\\
Bismarckstr. 10\\
76133 Karlsruhe\\
Germany\\
email: rosebrock@ph-karlsruhe.de

\end{document}